\documentclass[12pt]{article}

\usepackage{latexsym}
\usepackage{amssymb}
\usepackage{amsmath}
\usepackage{amsfonts}
\usepackage{txfonts}

 \numberwithin{equation}{subsection}

\begin{document}

\title{Characterizations of processes with stationary and independent increments under $G$-expectation }
\author{Yongsheng Song\\
\small Academy of Mathematics and Systems Science, \\
\small Chinese Academy of Sciences, Beijing, China;\\
\small yssong@amss.ac.cn}

\date{}

\maketitle
\begin{abstract}

Our purpose is to investigate  properties  for processes with
stationary and independent increments under $G$-expectation. As
applications, we prove the martingale  characterization  to
$G$-Brownian motion and present a decomposition theorem for
generalized $G$-Brownian motion.

\end{abstract}


{\bf MSC-classification}: 60G10, 60G17,  60G48, 60G51
\maketitle
\section{Introduction }

\ \ \  \ Recently, [P06], [P07], [P08a] introduced the notion of
$G$-expectation space, which is a generalization of probability
space. As the counterpart of Wiener space in the linear case, the
notions of $G$-Brownian motion, $G$-martingale, and It\^o integral
w.r.t. $G$-Brownian motion were also introduced.

 Recall that if $\{A_t\}$ is a continuous process in a probability space $(\Omega, {\cal F}, P)$ with stationary,
independent increments and finite variation, then there exists some
constant $c$ such that $A_t=ct$.  However, it's not the case in the
$G$-expectation space $(\Omega_T, L^1_G(\Omega_T), \hat{E})$. A
counterexample is $\{\langle B\rangle_t\}$, the quadratic variation
process for the coordinate process $\{B_t\}$, which is a
$G$-Brownian motion. We know that $\{\langle B\rangle_t\}$ is a
continuous, increasing process with stationary and independent
increments, however, it is not a deterministic process.

The process $\{\langle B\rangle_t\}$ is very important in the theory
of $G$-expectation, which shows, in many aspects, the difference
between probability space and $G$-expectation space. For example, we
know that in probability space continuous local martingales with
finite variation are trivial processes. However, [P07] proved that
in $G$-expectation space all processes in form of
$\int_0^t\eta_sd\langle B\rangle_s-\int_0^t2G(\eta_s)ds$, $\eta\in
M^1_G(0, T)$(see Section 2 for the definitions of the function
$G(\cdot)$ and the space $M^1_G(0, T)$), are $G$-martingales with
finite variation and continuous paths, which are a class of
nontrivial and very interesting processes. [P07] also conjectured
that any $G$-martingale with finite variation should have such
representation. Up to now, some properties of the process $\{\langle
B\rangle_t\}$ remain unknown. For example, we know that for any
$s<t$, $\underline{\sigma}^2(t-s)\leq\langle B\rangle_t-\langle
B\rangle_s\leq\overline{\sigma}(t-s)$, but we don't know whether
$\{\frac{d}{ds}\langle B\rangle_s\}$ belongs to $M^1_G(0, T)$. This
is a very important property since  $\{\frac{d}{ds}\langle
B\rangle_s\}\in M^1_G(0, T)$ implies that the representation
mentioned above of $G$-martingales with finite variation  is not
unique.

In probability space, a continuous local martingale $\{M_t\}$ with
the quadratic variation process $\langle M\rangle_t=t$ is a standard
Brownian motion. However, it's not the case for $G$-Brownian motion
since  its quadratic variation process is only an increasing process
with stationary and independent increments. How can we give a
characterization for $G$-Brownian motion?

In this article, we shall prove that if $A_t=\int_0^th_sds$
(respectively $A_t=\int_0^th_sd\langle B\rangle_s$) is a process
with stationary, independent increments and $h\in M^1_G(0,T)$
(respectively $h\in M^{1,+}_G(0,T)$), there exists some constant $c$
such that $h\equiv c$. As applications, we prove the following
conclusions (Question 1 and 3 are put forward by Prof. Shige Peng in
private communications):

1. $\{\frac{d}{ds}\langle B\rangle_s\}\notin M^1_G(0, T)$.

2.( Martingale characterization)

\textit{ A symmetric $G$-martingale $\{M_t\}$ is a $G$-Brownian
motion if and only if  its quadratic variation process $\{\langle
M\rangle_t\}$ has stationary and independent increments;}

\textit{A symmetric $G$-martingale $\{M_t\}$ is a $G$-Brownian
motion if and only if its quadratic variation process $\langle
M\rangle_t=c\langle B\rangle_t$ for some $c\geq0.$}

The sufficiency of the second part is implied by that of the first
part, but the necessity is not trivial.

3. \textit{Let $\{X_t\}$ be a generalized $G$-Brownian motion with zero
mean, then we have the following decomposition:
$$X_t=M_t+L_t,$$ where $\{M_t\}$ is a (symmetric) $G$-Brownian motion,
and $\{L_t\}$ is a non-positive, non-increasing $G$-martingale with
stationary and independent increments. }

This article is organized as follows: In section 2, we recall some
basic notions and results of $G$-expectation and the related space
of random variables. In section 3, we give characterizations to
processes with stationary and independent increments.  In section 4,
as applications, we prove the martingale  characterization to
$G$-Brownian motion and present a decomposition theorem for
generalized $G$-Brownian motion. In section 5, we present some
properties for $G$-martingales with finite variation.

\section{Preliminary }
We recall some basic notions and results of $G$-expectation and the
related space of random variables. More details of this section can
be found in [P06, P07, P08a, P08b, P10].

\noindent {\bf Definition 2.1} Let $\Omega$
 be a given set and let ${\cal H}$ be a vector lattice of real valued
functions defined on $\Omega$
 with $c \in {\cal H}$ for all constants $c$. ${\cal H}$ is considered as the
space of  ¡°random variables¡±. A sublinear expectation $\hat{E}$ on
${\cal H}$ is a functional $\hat{E}: {\cal H}\rightarrow R $
satisfying the following properties: For all $X, Y \in {\cal H}$, we
have

(a) Monotonicity: If $X\geq Y$ then $\hat{E}(X) \geq \hat{E} (Y)$.

(b) Constant preserving: $\hat{E} (c) = c$.

(c) Sub-additivity: $\hat{E}(X)-\hat{E}(Y) \leq \hat{E}(X-Y)$.

(d) Positive homogeneity: $\hat{E} (\lambda X) = \lambda \hat{E}
(X)$, $\lambda \geq 0$.

\noindent$(\Omega, {\cal H}, \hat{E})$ is called a sublinear
expectation space. $\Box$

\noindent {\bf Definition 2.2} Let $X_1$ and $X_2$ be two
$n$-dimensional random vectors defined respectively in sublinear
expectation spaces $(\Omega_1, {\cal H}_1, \hat{E}_1)$ and
$(\Omega_2, {\cal H}_2, \hat{E}_2)$. They are called identically
distributed, denoted by $X_1 \sim X_2$, if $\hat{E}_1[\varphi(X_1)]
= \hat{E}_2[\varphi(X_2)]$, $\textmd{for all} \ \varphi\in C_{l,
Lip}(R^n)$, where $ C_{l, Lip}(R^n)$ is the space of real continuous
functions defined on $R^n$ such that $$|\varphi(x) - \varphi(y)|
\leq C(1 + |x|^k + |y|^k)|x - y|, \textmd{for all} \ x, y \in R^n,$$
where $k$ and $C$ depend only on $\varphi$. $\Box$

\noindent {\bf Definition 2.3} In a sublinear expectation space
$(\Omega, {\cal H}, \hat{E})$ a random vector $Y = (Y_1,
\cdot\cdot\cdot, Y_n)$, $Y_i \in {\cal H}$, is said to be
independent of another random vector $X = (X_1, \cdot\cdot\cdot,
X_m)$, $X_i \in {\cal H}$, under $\hat{E}(\cdot)$, denoted by
$Y\Perp X$, if for every test function $\varphi\in C_{b,
Lip}(R^m\times R^n)$ we have $\hat{E}[\varphi(X, Y )] =
\hat{E}[\hat{E} [\varphi(x, Y )]_{x=X}]$. $\Box$

\noindent {\bf Definition 2.4} ($G$-normal distribution) A
d-dimensional random vector $X = (X_1, \cdot\cdot\cdot,X_d)$ in a
sublinear expectation space $(\Omega, {\cal H}, \hat{E})$ is called
$G$-normal distributed if for every $a, b\in R_+$ we have $$aX +
b\hat{X}\sim \sqrt{a^2 + b^2}X,$$  where $\hat{X}$ is an independent
copy of $X$. Here the letter $G$ denotes the function $$G(A) :=
\frac{1 }{2}\hat{ E}[(AX,X)] : S_d \rightarrow R,$$  where $S_d$
denotes the collection of $d\times d$ symmetric matrices. $\Box$

The function $G(\cdot) : S_d \rightarrow R$ is a monotonic,
sublinear mapping on $S_d$ and $G(A) = \frac{1 }{2}\hat{
E}[(AX,X)]\leq \frac{1 }{2}|A|\hat{ E}[|X|^2]=:\frac{1
}{2}|A|\bar{\sigma}^2$ implies that there exists a bounded, convex
and closed subset $\Gamma\subset S_d^+$ such that
\begin {eqnarray}\label{eqn1}
G(A)=\frac{1 }{2}\sup_{\gamma\in \Gamma}Tr(\gamma A).
\end {eqnarray}
If there exists some $\beta>0$ such that $G(A)-G(B)\geq \beta
Tr(A-B)$ for any $A\geq B$, we call the $G$-normal distribution
non-degenerate. This is the case we consider throughout this
article.

\noindent {\bf Definition 2.5} i) Let $\Omega_T=C_0([0, T]; R^d)$ be
endowed with the supremum norm and $\{B_t\}$ be the coordinate
process.  Set $ {\cal H}^0_T:=\{\varphi(B_{t_1},..., B_{t_n})|
n\geq1, t_1, ..., t_n \in [0, T],  \varphi \in C_{l, Lip}(R^{d\times
n})\}$.  $G$-expectation is a sublinear expectation defined by
$$\hat{E}[X] = \tilde{E}
[\varphi(\sqrt{t_1-t_0}\xi_1, \cdot\cdot\cdot, \sqrt{t_m
-t_{m-1}}\xi_m)],$$ for all $X=\varphi( B_{t_1}-B_{t_0},
B_{t_2}-B_{t_1} , \cdot\cdot\cdot, B_{t_m}- B_{t_{m-1}} )$, where
$\xi_1, \cdot\cdot\cdot, \xi_n$ are identically distributed
$d$-dimensional $G$-normally distributed random vectors in a
sublinear expectation space $(\tilde{\Omega}, \tilde{\cal H},\tilde{
E})$ such that  $\xi_{i+1}$ is independent of $(\xi_1,
\cdot\cdot\cdot, \xi_i)$ for every $i = 1, \cdot\cdot\cdot,m$.
$(\Omega_T, {\cal H}^0_T, \hat{E})$ is called a $G$-expectation
space.

ii) For $t\in [0, T]$ and $\xi=\varphi(B_{t_1},..., B_{t_n})\in
{\cal H}^0_T$, the conditional expectation defined by(there is no
loss of generality, we assume $t=t_i$) $$\hat{E}_{t_i}[\varphi(
B_{t_1}-B_{t_0}, B_{t_2}-B_{t_1} , \cdot\cdot\cdot, B_{t_m}-
B_{t_{m-1}} )]$$$$=\tilde{\varphi}( B_{t_1}-B_{t_0}, B_{t_2}-B_{t_1}
, \cdot\cdot\cdot, B_{t_i}- B_{t_{i-1}} ),$$ where
$$\tilde{\varphi}(x_1, \cdot\cdot\cdot, x_i)=\hat{E}[\varphi( x_1,
\cdot\cdot\cdot,x_i, B_{t_{i+1}}- B_{t_{i}}, \cdot\cdot\cdot,
B_{t_m}- B_{t_{m-1}} )].$$ $\Box$

Define $\|\xi\|_{p, G}=[\hat{E}(|\xi|^p)]^{1/p}$ for $\xi\in{\cal
H}^0_T$ and $p\geq1$. Then for all $t\in[0, T]$, $\hat{E}_t(\cdot)$
is a continuous mapping on ${\cal H}^0_T$ with respect to the norm
$\|\cdot\|_{1, G}$ and therefore can be extended continuously to the
completion $L^1_G(\Omega_T)$ of ${\cal H}^0_T$ under norm
$\|\cdot\|_{1, G}$.

Let $L_{ip}(\Omega_T):=\{\varphi(B_{t_1},..., B_{t_n})|  n\geq1,
t_1, ..., t_n \in [0, T],  \varphi \in C_{b, Lip}(R^{d\times n})\},$
where $C_{b, Lip}(R^{d\times n})$ denotes the set of bounded
Lipschitz functions on $R^{d\times n}$. [DHP08] proved that the
completions of $C_b(\Omega_T)$, ${\cal H}^0_T$ and
$L_{ip}(\Omega_T)$ under $\|\cdot\|_{p,G}$ are  same and we denote
them by $L^p_G(\Omega_T)$.

\noindent {\bf Definition 2.6} i) We  say  that $\{X_t\}$ on
$(\Omega_T, L^1_G(\Omega_T), \hat{E})$ is a process with independent
increments if for any $0<t<T$ and $s_0\leq\cdot\cdot\cdot\leq
s_m\leq t\leq t_0\leq\cdot\cdot\cdot\leq t_n\leq T$,
$$(X_{t_1}-X_{t_{0}}, \cdot\cdot\cdot,
X_{t_n}-X_{t_{n-1}})\Perp (X_{s_1}-X_{s_{0}}, \cdot\cdot\cdot,
X_{s_m}-X_{s_{m-1}}).$$

ii) We say  that $\{X_t\}$ on $(\Omega_T, L^1_G(\Omega_T), \hat{E})$
with $X_t\in L^1_G(\Omega_t)$ for every $t\in[0,T]$ is a process
with independent increments {\bf w.r.t. the filtration } if for any
$0<s<T$ and $s_0\leq\cdot\cdot\cdot\leq s_m\leq s\leq
t_0\leq\cdot\cdot\cdot \leq t_n\leq T$, $$(X_{t_1}-X_{t_{0}},
\cdot\cdot\cdot, X_{t_n}-X_{t_{n-1}}) \Perp (B_{s_1}-B_{s_{0}},
\cdot\cdot\cdot, B_{s_m}-B_{s_{m-1}}).$$ $\Box$

\noindent {\bf Remark 2.7} i)  Let $\xi\in L^1_G(\Omega_T)$. If
there exists $s\in[0,T]$ such that for any
$s_0\leq\cdot\cdot\cdot\leq s_m\leq s$, $\xi\Perp
(B_{s_1}-B_{s_{0}}, \cdot\cdot\cdot, B_{s_m}-B_{s_{m-1}}),$ then we
have $\hat{E}_s(\xi)=\hat{E}(\xi)$. In fact, there is no loss of
generality, we assume $\xi\geq0$.
\begin {eqnarray*}& &\hat{E}[\hat{E}_s(\xi)-\hat{E}(\xi)]^2\\
&=&\hat{E}[\hat{E}_s(\xi)(\hat{E}_s(\xi)-2\hat{E}(\xi))]+[\hat{E}(\xi)]^2\\
&=&\hat{E}[\hat{E}_s(\xi)(\xi-2\hat{E}(\xi))]+[\hat{E}(\xi)]^2.
\end {eqnarray*} Since $\xi-2\hat{E}(\xi)\Perp\hat{E}_s(\xi)$, we
have
$$\hat{E}[\hat{E}_s(\xi)-\hat{E}(\xi)]^2=\hat{E}(\xi)(\hat{E}(\xi)-2\hat{E}(\xi))+[\hat{E}(\xi)]^2=0.$$
ii) Let $\{X_t\}$ on $(\Omega_T, L^1_G(\Omega_T), \hat{E})$ be a
process with stationary and independent increments and let
$c=\hat{E}(X_T)/T$. If $\hat{E}(X_t)\rightarrow0$ as $t\downarrow0$,
then for any $0\leq s<t\leq T$, we have
$\hat{E}(X_t-X_s)=c(t-s)$.$\Box$

\noindent {\bf Definition 2.8} Let $\{X_t\}$ be a d-dimensional
process defined on $(\Omega_T, L^1_G(\Omega_T), \hat{E})$ such that

(i) $X_0 = 0$;

(ii) $\{X_t\}$ is a process with stationary and independent
increments w.r.t. the filtration;

(iii) $\lim_{t\rightarrow0} \hat{E}[|X_t|^3]t^{-1} = 0.$

Then $\{X_t\}$ is called a generalized $G$-Brownian motion.

If in addition $\hat{E}(X_t)=\hat{E}(-X_t)=0$, $\{X_t\}$ is called a
(symmetric) $G$-Brownian motion. $\Box$

\noindent {\bf Remark 2.9} i) Clearly, the coordinate process
$\{B_t\}$ is a (symmetric) $G$-Brownian motion and its quadratic
variation process $\{\langle B\rangle_t\}$ is a process with
stationary and independent increments (w.r.t. the filtration).

ii) [P07] gave a characterization for the generalized $G$-Brownian
motion: Let $\{X_t\}$ be a generalized $G$-Brownian motion. Then
\begin {eqnarray} \label {eqns7}X_{t+s}-X_t\sim \sqrt{s}\xi + s\eta \ \textmd{for} \ t,
s\geq0, \end {eqnarray} where $(\xi, \eta)$ is $G$-distributed(see,
e.g., [P08b] for the definition of $G$-distributed random vectors).
In fact, the characterization presented a decomposition of
generalized $G$-Brownian motion in the sense of distribution. In
this article, we shall give a pathwise decomposition for  the
generalized $G$-Brownian motion. $\Box$

Let $H^0_G(0, T)$ be the collection of processes of the following
form: for a given partition $\{t_0, \cdot\cdot\cdot, t_N\} = \pi_T$
of $[0, T]$, $N\geq1$,  $$ \eta_t(\omega) = \sum^{N-1}_{j=0}
\xi_j(\omega)1_{]t_j ,t_{j+1}]}(t),$$ where $\xi_i\in L_{ip}(\Omega_
{t_i})$, $i = 0, 1, 2, \cdot\cdot\cdot, N-1$. For every $\eta\in
H^0_G(0, T)$, let
$\|\eta\|_{H^{p}_G}=\{\hat{E}(\int_0^T|\eta_s|^2ds)^{p/2}\}^{1/p}$,
$\|\eta\|_{M^{p}_G}=\{\hat{E}(\int_0^T|\eta_s|^pds)\}^{1/p}$ and
denote by $H^{p}_G(0, T)$,  $M^{p}_G(0, T)$ the completions of
$H^0_G(0, T)$ under the norms $\|\cdot\|_{H^{p}_G}$,
$\|\cdot\|_{M^{p}_G}$ respectively.

\noindent {\bf Definition 2.10} For every $\eta\in H^0_G(0, T)$ with
the form $$\eta_t(\omega) = \sum^{N-1}_{j=0}
\xi_j(\omega)1_{]t_j,t_{j+1}]}(t),$$ we define $$I(\eta) =\int_0^T
\eta(s)dB_s := \sum^{N-1}_{j=0} \xi_j(B_{t_{j+1} }-B_{t_j} ).$$

By B-D-G inequality(see Proposition 4.3 in [song11a] for this
inequality under $G$-expectation), the mapping $I: H^0_G(0,
T)\rightarrow L^p_G(\Omega_T)$ is continuous under
$\|\cdot\|_{H^{p}_G}$ and thus can be continuously extended to
$H^p_G(0, T)$. $\Box$

\noindent {\bf Definition 2.11} i) A process $\{M_t\}$ with values
in $L^1_G(\Omega_T)$ is called a $G$-martingale if
$\hat{E}_s(M_t)=M_s$ for any $s\leq t$. If $\{M_t\}$ and  $\{-M_t\}$
are both $G$-martingales, we call $\{M_t\}$ a symmetric
$G$-martingale.

ii) A random variable $\xi\in L^1_G(\Omega_T)$ is called symmetric,
if $\hat{E}(\xi)+\hat{E}(-\xi)=0.$ $\Box$

A $G$-martingale $\{M_t\}$ is symmetric  if and only if $M_T$ is
symmetric.

\noindent {\bf Theorem 2.12}([DHP08]) There exists a tight subset
${\cal P}\subset {\cal M}_1(\Omega_T)$ such that
$$\hat{E}(\xi)=\max_{P\in {\cal P}}E_P(\xi) \ \ \textrm{for \
all} \ \xi\in{\cal H}^0_T.$$ ${\cal P}$ is called a set that
represents $\hat{E}$.

\noindent {\bf Remark 2.13} i) Let $(\Omega^0, {\cal F}^0, P^0 )$ be
a
 probability space and $\{W_t\}$ be a d-dimensional Brownian
motion under $P^0$. Let $F^0=\{{\cal F}^0_t\}$ be the augmented
filtration generated by $W$.  [DHP08] proved that
$${\cal P}_M:=\{P_h|P_h=P^0\circ X^{-1}, X_t=\int_0^th_sdW_s, h\in
L^2_{ F^0}([0,T]; \Gamma^{1/2}) \}$$ is a set that represents
$\hat{E}$, where $\Gamma^{1/2}:=\{\gamma^{1/2}| \gamma\in \Gamma\}$
and $\Gamma$ is the set in the representation of $G(\cdot)$ in the
formula (\ref {eqn1}).

ii) For the 1-dimensional case, i.e., $\Omega_T=C_0([0,T], R^1)$,
$$L^2_{F^0}:=L^2_{ F^0}([0,T]; \Gamma^{1/2})=\{h| \ h \ \textmd{is
adapted w.r.t.} \  F^0 \ \textmd{and} \ \underline{\sigma}\leq
h_s\leq\overline{\sigma}\},$$ where
$\overline{\sigma}^2=\hat{E}(B_1^2)$ and $\underline{\sigma}^2=
-\hat{E}(-B_1^2)$.

$G(a)=1/2\hat{E}[aB_1^2]=1/2[\overline{\sigma}^2a^+-\underline{\sigma}^2a^-]$
for $a\in R.$

iii) Set $c(A)=\sup_{P\in{\cal P}_M}P(A)$, for $A\in {\cal
B}(\Omega_T)$. We say $A\in {\cal B}(\Omega_T)$ is a polar set if
$c(A)=0$. If an event happens except on a polar set, we say the
event happens q.s..

\section{Characterization of processes with stationary and independent increments }
In what follows, we only consider the $G$-expectation space
$(\Omega_T, L^1_G(\Omega_T), \hat{E})$ with $\Omega_T=C_0([0,T], R)$
and
$\overline{\sigma}^2=\hat{E}(B_1^2)>-\hat{E}(-B_1^2)=\underline{\sigma}^2>0$.

\noindent {\bf Lemma 3.1} For $\zeta\in M^1_G(0,T)$ and
$\varepsilon>0$, let
$$\zeta^\varepsilon_t=\frac{1}{\varepsilon}\int^t_{(t-\varepsilon)^+}\zeta_sds$$
and
$$\zeta^{\varepsilon,0}_t=\sum_{k=1}^{k_\varepsilon-1}\frac{1}{\varepsilon}\int^{k\varepsilon}_{(k-1)\varepsilon}\zeta_sds1_{]k\varepsilon,
(k+1)\varepsilon]}(t),$$ where $t\in[0,T]$,
$k_\varepsilon\varepsilon\leq T<(k_\varepsilon+1)\varepsilon$. Then
as $\varepsilon\rightarrow0$
$$\|\zeta^\varepsilon-\zeta\|_{M^1_G(0,T)}\rightarrow0 \ \ \textmd{and} \ \
\|\zeta^{\varepsilon,0}-\zeta\|_{M^1_G(0,T)}\rightarrow0.$$

{\bf Proof.} The proofs of the two cases are similar. Here we only
prove the second case. Our proof starts with the observation that
for any $\zeta, \zeta'\in M^1_G(0,T)$
\begin{eqnarray}\label{eqns2}
\|\zeta^{\varepsilon,0}-\zeta'^{\varepsilon,0}\|_{M^1_G(0,T)}\leq\|\zeta-\zeta'\|_{M^1_G(0,T)}.
\end {eqnarray}
By the definition of space $M^1_G(0,T)$, we  know that  for every
$\zeta\in M^1_G(0,T)$, there exists a sequence of processes
$\{\zeta^n\}$ with
$$\zeta^n_t=\Sigma_{k=0}^{m_n-1}\xi^n_{t^n_k}1_{]t^n_k, t^n_{k+1}]}(t)$$
 and $\xi^n_{t^n_k}\in L_{ip}(\Omega_{t^n_k})$ such that
\begin{eqnarray}\label{eqns3}
\|\zeta-\zeta^n\|_{M^1_G(0,T)}\rightarrow0 \ \textmd{as} \
n\rightarrow\infty.
\end {eqnarray}
 It is easily seen that
for every $n$,
\begin {eqnarray}\label{eqns4} \|\zeta^{n;\varepsilon,0}-\zeta^n\|_{M^1_G(0,T)}\rightarrow0
\ \textmd{as} \ \varepsilon\rightarrow0.
\end {eqnarray}
Thus we get
\begin{eqnarray*}& &\|\zeta^{\varepsilon,0}-\zeta\|_{M^1_G(0,T)}\\
                 &\leq&\|\zeta^{\varepsilon,0}-\zeta^{n;\varepsilon,0}\|_{M^1_G(0,T)}
+\|\zeta^n-\zeta^{n;\varepsilon,0}\|_{M^1_G(0,T)}+\|\zeta^n-\zeta\|_{M^1_G(0,T)}\\
                 &\leq&2\|\zeta^n-\zeta\|_{M^1_G(0,T)}+\|\zeta^n-\zeta^{n;\varepsilon,0}\|_{M^1_G(0,T)}.
\end {eqnarray*}
The second inequality follows from (\ref {eqns2}). Combining (\ref
{eqns3}) and (\ref {eqns4}), first letting
$\varepsilon\rightarrow0$, then letting $n\rightarrow\infty$, we
have
$$\|\zeta^{\varepsilon,0}-\zeta\|_{M^1_G(0,T)}\rightarrow0 \ \textmd{as} \ \varepsilon\rightarrow0.$$ $\Box$

\noindent {\bf Theorem 3.2} Let $A_t=\int_0^th_sds$ with $h\in
M^1_G(0, T)$ be a process with stationary and independent increments
(w.r.t. the filtration). Then we have $h\equiv c$ for some constant
$c$.

{\bf Proof.} Let
$\overline{c}:=\hat{E}(A_T)/T\geq-\hat{E}(-A_T)/T=:\underline{c}$.
For  $n\in N$, set $\varepsilon=T/(2n)$, and define $h^{T/(2n),0}$
as in Lemma 3.1. Then we have
\begin{eqnarray*}& &\|h-h^{T/(2n),0}\|_{M^1_G(0,T)}\\
                 &=&\hat{E}[\sum_{k=0}^{2n-1}\int_{kT/(2n)}^{(k+1)T/(2n)}|h_s-h^{T/(2n),0}_s|ds]\\
                 &\geq&\hat{E}[\sum_{k=1}^{n-1}\int_{2kT/(2n)}^{(2k+1)T/(2n)}(h_s-h^{T/(2n),0}_s)ds]\\
                 &=&\hat{E}[\sum_{k=1}^{n-1}(\int_{2kT/(2n)}^{(2k+1)T/(2n)}h_sds-\int_{(2k-1)T/(2n)}^{2kT/(2n)}h_sds)]\\
                 &=&\hat{E}{\sum_{k=1}^{n-1}[(
                 A_{(2k+1)T/2n}-A_{2kT/2n})-(A_{2kT/2n}-A_{(2k-1)T/2n})]}.
\end{eqnarray*}
Consequently, from the condition of independence of the increments
and their stationarity,

\begin{eqnarray*}& &\|h-h^{T/(2n),0}\|_{M^1_G(0,T)}\\
                 &\geq&\sum_{k=1}^{n-1}\hat{E}[(
                 A_{(2k+1)T/2n}-A_{2kT/2n})-(A_{2kT/2n}-A_{(2k-1)T/2n})]\\
                 &=&\sum_{k=1}^{n-1}(\overline{c}-\underline{c})T/(2n)\\
                 &=&(\overline{c}-\underline{c})(n-1)T/(2n).
\end{eqnarray*}

So by Lemma 3.1, letting $n\rightarrow\infty$, we have
$\overline{c}=\underline{c}$. Furthermore, we note that
$M_t:=A_t-\overline{c}t$ is a $G$-martingale. In fact, for $t>s$, we
see
\begin {eqnarray*}& &\hat{E}_s(M_t)\\
                 &=&\hat{E}_s(M_t-M_s)+M_s\\
                 &=&\hat{E}(M_t-M_s)+M_s\\
                 &=&M_s.
\end {eqnarray*} The second equality is due to the independence of increments of $M$ w.r.t. the filtration.

So $\{M_t\}$ is a symmetric $G$-martingale with finite variation,
from which we conclude that $M_t\equiv0$, hence that
$A_t=\overline{c}t$.$\Box$

\noindent {\bf Corollary 3.3} Assume
$\overline{\sigma}>\underline{\sigma}>0$. Then we have that
$\{\frac{d}{ds}\langle B\rangle_s\}\notin M^1_G(0, T)$.

{\bf Proof.} The proof is straightforward from Theorem 3.2. $\Box$

\noindent {\bf Corollary 3.4} There is no symmetric $G$-martingale
$\{M_t\}$ which is  a standard Brownian motion under
$G$-expectation(i.e. $\langle M\rangle_t=t$).

{\bf Proof.} Let $\{M_t\}$ be a symmetric $G$-martingale. If
$\{M_t\}$ is also a standard Brownian motion, by Theorem 4.8 in
[Song11a] or Corollary 5.2 in [Song11b], there exists $\{h_s\}\in
M^2_G(0, T)$ such that
$$M_t=\int_0^th_sdB_s$$ and $$\int_0^th_s^2d\langle B\rangle_s=t.$$
Thus we have $\frac{d}{ds}\langle B\rangle_s=h_s^{-2}\in
M^1_G(0,T)$, which contradicts the conclusion of Corollary 3.3.
$\Box$

\noindent {\bf Proposition 3.5} Let $A_t=\int_0^th_sds$ with $h\in
M^1_G(0, T)$ be a process with  independent increments. Then $A_t$
is symmetric for every $t\in[0,T]$.

{\bf Proof.}  By  arguments similar to  that in the proof of Theorem
3.2, we have

\begin{eqnarray*}& &\|h-h^{T/(2n),0}\|_{M^1_G(0,T)}\\
                 &\geq&\hat{E}{\sum_{k=0}^{n-1}[(
                 A_{(2k+1)T/2n}-A_{2kT/2n})-(A_{2kT/2n}-A_{(2k-1)^+T/2n})]}\\
                 &=&\sum_{k=0}^{n-1}\{\hat{E}(
                 A_{(2k+1)T/2n}-A_{2kT/2n})+\hat{E}[-(A_{2kT/2n}-A_{(2k-1)^+T/2n})]\}.
\end{eqnarray*}The right side of the first inequality is only the sum of the odd
terms.  Summing up  the even terms only, we have

\begin{eqnarray*}& &\|h-h^{T/(2n),0}\|_{M^1_G(0,T)}\\
                 &\geq&\sum_{k=0}^{n-1}\{\hat{E}(
                 A_{(2k+2)T/2n}-A_{(2k+1)T/2n})+\hat{E}[-(A_{(2k+1)T/2n}-A_{2kT/2n})]\}.
\end{eqnarray*}
Combining the above inequalities, we have

\begin{eqnarray*}& &2\|h-h^{T/(2n),0}\|_{M^1_G(0,T)}\\
                 &\geq&\sum_{k=0}^{2n-1}\{\hat{E}(
                 A_{(k+1)T/2n}-A_{(kT/2n})+\hat{E}[-(A_{(k+1)T/2n}-A_{kT/2n})]\}\\
                 &\geq&\hat{E}\sum_{k=0}^{2n-1}(
                 A_{(k+1)T/2n}-A_{(kT/2n})+\hat{E}\sum_{k=0}^{2n-1}[-(A_{(k+1)T/2n}-A_{kT/2n})]\\
                 &=& \hat{E}(A_T)+\hat{E}(-A_T).
\end{eqnarray*}

Thus by Lemma 3.1, letting $n\rightarrow\infty$, we have
$\hat{E}(A_T)+\hat{E}(-A_T)=0$, which means that $A_T$ is symmetric.
$\Box$

For $n\in N$, define $\delta_n(s)$ in the following way:
$$\delta_n(s)=\sum_{i=0}^{n-1}(-1)^i1_{]\frac{iT}{n}, \frac{(i+1)T}{n}]}(s), \textmd{\ for \ all} \
s\in[0,T].$$ In [Song10c] we proved that
$\lim_{n\rightarrow\infty}\hat{E}(\int_0^T\delta_n(s)h_sds)=0$ for
$h\in M^1_G(0,T)$.

Let ${\cal F}_t=\sigma\{B_s| s\leq t\}$ and $\mathbb{F}=\{{\cal
F}_t\}_{t\in[0, T]}$.

In the following, we shall use some notations introduced in Remark
2.13.

For every $P\in{\cal P}_M$ and $t\in[0, T]$, set ${\cal A}_{t,
P}:=\{Q\in{\cal P}_M| \ Q_{|{{\cal F}_t}}=P_{|{{\cal F}_t}}\}$.
Proposition 3.4 in [STZ11] gave the following result: For $t\in[0,
T]$, assume $\xi\in L^1_G(\Omega_T)$ and $\eta\in L^1_G(\Omega_t)$.
Then $\eta=\hat{E}_t(\xi)$ if and only if for every $P\in{\cal P}_M$
$$\eta={\textmd{ess}\sup}_{Q\in {\cal A}_{t,
P}}^PE_Q(\xi|{\cal F}_t), \ \textit{P-a.s.},$$ where
${\textmd{ess}\sup}^P$ denotes the essential supremum under $P$.

\noindent {\bf Theorem 3.6} Let $A_t=\int_0^th_sd\langle B\rangle_s$
be a process with stationary, independent increments (w.r.t. the
filtration) and $h\in M^{\beta,+}_G(0,T)$ for some $\beta>1$. Then
there exists a constant $c\geq0$ such that $A_t=c\langle
B\rangle_t$.

{\bf Proof.} For the readability, we divide the proof into several
steps:

{\bf Step 1.} Set $K_t:=\int_0^th_sds$. We claim that $K_T$ is
symmetric.

{\bf Step 1.1.} Let $\overline{\mu}=\hat{E}(A_T)/T$ and
$\underline{\mu}=-\hat{E}(-A_T)/T$. First, we shall prove that
$\frac{\overline{\mu}}{\overline{\sigma}^2}=\frac{\underline{\mu}}{\underline{\sigma}^2}$.

Actually, for any $0\leq s<t\leq T$, we have
$$\hat{E}_s(\int_s^th_rdr)=\hat{E}_s(\int_s^t\theta_r^{-1}dA_r)\geq\frac{1}{\overline{\sigma}^2}\hat{E}_s(\int_s^tdA_r)=
\frac{\overline{\mu}}{\overline{\sigma}^2}(t-s)\ \ \ q.s.,$$ where
the inequality holds due to $\theta_s\leq\overline{\sigma}^2$, q.s..
Since $h\in M^\beta_G(\Omega_T)$, we have $A_T\in
L^\beta_G(\Omega_T)$. Noting that $\underline{\mu}t-A_t$ is
nonincreasing by Lemma 4.3 in Section 4 since it is a $G$-martingale
with finite variation, we have, for every $\eta\in L^2_{F^0}$,
$P_\eta$-a.s.,
\begin {eqnarray*}
& &\hat{E}_s(\int_s^th_rdr)\\
&=&{\textmd{ess}\sup}_{Q\in {\cal A}_{t,P_\eta}}^{P_\eta}E_Q(\int_s^th_rdr|{\cal F}_s)\\
&=&{\textmd{ess}\sup}_{Q\in {\cal
A}_{t,P_\eta}}^{P_\eta}E_Q(\int_s^t\theta_r^{-1}dA_r|{\cal
F}_s)\\
&\geq&\underline{\mu} \ {\textmd{ess}\sup}_{Q\in {\cal
A}_{t,P_\eta}}^{P_\eta}E_Q(\int_s^t\theta_r^{-1}dr|{\cal F}_s)\\
&=&\frac{\underline{\mu}}{\underline{\sigma}^2}(t-s).
\end {eqnarray*}

 So
$$\hat{E}_s(\int_s^th_rdr)\geq \max\{\frac{\overline{\mu}}{\overline{\sigma}^2},\frac{\underline{\mu}}{\underline{\sigma}^2}\}(t-s)=:\overline{\lambda}(t-s), \ q.s..$$
On the other hand,
$$\hat{E}_s(-\int_s^th_rdr)=\hat{E}_s(\int_s^t-\theta_r^{-1}dA_r)\geq\frac{1}{\underline{\sigma}^2}\hat{E}_s(-\int_s^tdA_r)=-
\frac{\underline{\mu}}{\underline{\sigma}^2}(t-s), \ q.s.$$ and for
every $\eta\in L^2_{F^0}$, $P_\eta$-a.s.,
\begin {eqnarray*}
& &\hat{E}_s(-\int_s^th_rdr)\\
&=&{\textmd{ess}\sup}_{Q\in {\cal A}_{t,P_\eta}}^{P_\eta}E_Q(-\int_s^th_rdr|{\cal F}_s)\\
&=&{\textmd{ess}\sup}_{Q\in {\cal
A}_{t,P_\eta}}^{P_\eta}E_Q(-\int_s^t\theta_r^{-1}dA_r|{\cal
F}_s)\\
&\geq&\overline{\mu} \ {\textmd{ess}\sup}_{Q\in {\cal
A}_{t,P_\eta}}^{P_\eta}E_Q(-\int_s^t\theta_r^{-1}dr|{\cal F}_s)\\
&=&-\frac{\overline{\mu}}{\overline{\sigma}^2}(t-s)
\end {eqnarray*} since $A_t-\overline{\mu}t$ is nonincreasing. So
$$\hat{E}_s(-\int_s^th_rdr)\geq -\min\{\frac{\overline{\mu}}{\overline{\sigma}^2},\frac{\underline{\mu}}{\underline{\sigma}^2}\}(t-s)=:-\underline{\lambda}(t-s), \ q.s..$$
Noting that
\begin {eqnarray*}
& & \hat{E}(\int_0^T\delta_{2n}(s)h_sds)\\
&=&
\hat{E}[\int_0^{\frac{(2n-1)T}{2n}}\delta_{2n}(s)h_sds+\hat{E}_{\frac{(2n-1)T}{2n}}(-\int_{\frac{(2n-1)T}{2n}}^Th_sds)]\\
&\geq&(-\underline{\lambda})\frac{T}{2n}+\hat{E}[\int_0^{\frac{(2n-2)T}{2n}}\delta_{2n}(s)h_sds+
\hat{E}_{\frac{(2n-2)T}{2n}}(\int_{\frac{(2n-2)T}{2n}}^{\frac{(2n-1)T}{2n}}h_sds)]\\
&\geq&\frac{\overline{\lambda}-\underline{\lambda}}{2n}T+\hat{E}[\int_0^{\frac{(2n-2)T}{2n}}\delta_{2n}(s)h_sds],
\end {eqnarray*} we have
$$\hat{E}(\int_0^T\delta_{2n}(s)h_sds)\geq\frac{\overline{\lambda}-\underline{\lambda}}{2}T.$$
So
$$0=\lim_{n\rightarrow\infty}\hat{E}(\int_0^T\delta_{2n}(s)h_sds)\geq\frac{\overline{\lambda}-\underline{\lambda}}{2}T,$$
and
$\frac{\overline{\mu}}{\overline{\sigma}^2}=\frac{\underline{\mu}}{\underline{\sigma}^2}=:\lambda$.

{\bf Step 1.2.} For every $\eta\in L^2_{F^0}$,
$E_{P_\eta}(K_T)=\lambda T$, which implies that $K_T$ is symmetric.

{\bf Step 1.2.1.} We now introduce some notations: For $0\leq
s<t\leq T$ and $\eta\in L^2_{F^0}$, set
$\overline{\eta}=\overline{\sigma}$,
$\underline{\eta}=\underline{\sigma}$,
$\eta^*=\sqrt{\frac{\overline{\sigma}^2+\underline{\sigma}^2}{2}}$
on $]s,t]$ and $\overline{\eta}=\underline{\eta}=\eta^*=\eta$ on
$]s,t]^c$. For $n\in N$, set
$\eta^n_r=\sum_{i=0}^{n-1}(\underline{\sigma}1_{]t_{2i},
t_{2i+1}]}(r)+\overline{\sigma}1_{]t_{2i+1}, t_{2i+2}]}(r))$ on
$]s,t]$ and $\eta^n=\eta$ on $]s,t]^c$, where
$t_j=s+\frac{j}{2n}(t-s)$, $j=0, \cdot\cdot\cdot, 2n$.

{\bf Step 1.2.2.} $E_{P_{\eta^n}}(\int_s^t(h_r-\lambda)dr|{\cal
F}_s)\rightarrow0$, $P_\eta$-a.s., as $n\rightarrow\infty$.

Actually, we have, $P_\eta$-a.s.,
$$\overline{\mu}(t-s)=\hat{E}_s(\int_s^th_rd\langle B\rangle_r)\geq E_{P_{\overline{\eta}}}(\int_s^th_rd\langle B\rangle_r|{\cal F}_s)
=\overline{\sigma}^2E_{P_{\overline{\eta}}}(\int_s^th_rdr|{\cal
F}_s).$$ So
\begin {eqnarray}\label{eqns5}E_{P_{\overline{\eta}}}(\int_s^th_rdr|{\cal F}_s)\leq
\lambda(t-s), \ P_\eta-a.s..
\end {eqnarray}

 By similar arguments we have that
\begin {eqnarray}\label{eqns6}E_{P_{\underline{\eta}}}(\int_s^th_rdr|{\cal F}_s)\geq
\lambda(t-s), \ P_\eta-a.s..
\end {eqnarray}

Let's compute the following conditional expectations:
\begin {eqnarray*}
& &E_{P_{\eta^n}}(\int_s^t(h_r-\lambda)\delta_{2n}(r)dr|{\cal F}_s)\\
&=&E_{P_{\eta^n}}^{{\cal
F}_s}[\sum_{i=0}^{n-1}\{E_{P_{\eta^n}}^{{\cal
F}_{t_{2i}}}\int_{t_{2i}}^{t_{2i+1}}(h_r-\lambda)dr+E_{P_{\eta^n}}^{{\cal
F}_{t_{2i+1}}}\int_{t_{2i+1}}^{t_{2i+2}}(\lambda-h_r)dr\}]\\
&=&: E_{P_{\eta^n}}^{{\cal F}_s}[\sum_{i=0}^{n-1}(A_i+B_i)],
\end{eqnarray*} where $\delta_{2n}(r)=\sum_{i=0}^{n-1}(1_{]t_{2i},
t_{2i+1}]}(r)-1_{]t_{2i+1}, t_{2i+2}]}(r))$,
$t_j=s+\frac{j}{2n}(t-s)$, $j=0, \cdot\cdot\cdot, 2n$, and
\begin {eqnarray*}
E_{P_{\eta^n}}(\int_s^t(h_r-\lambda)dr|{\cal F}_s)=
E_{P_{\eta^n}}^{{\cal F}_s}[\sum_{i=0}^{n-1}(A_i-B_i)].
\end{eqnarray*}

By (\ref {eqns5}) and (\ref {eqns6})(noting that $\eta$ and $s,t$
are all arbitrary), we conclude that $A_i, B_i\geq0$,
$P_{\eta_n}$-a.s.. So
$$|E_{P_{\eta^n}}(\int_s^t(h_r-\lambda)dr|{\cal F}_s)|\leq E_{P_{\eta^n}}(\int_s^t(h_r-\lambda)\delta_{2n}(r)dr|{\cal F}_s), \ P_\eta-a.s..$$

Noting that
$$E_{P_{\eta^n}}(\int_s^t(h_r-\lambda)\delta_{2n}(r)dr|{\cal
F}_s)\leq \hat{E}_s[\int_s^t(h_r-\lambda)\delta_{2n}(r)dr], \
P_\eta-a.s.$$ and
$$\hat{E}_s[\int_s^t(h_r-\lambda)\delta_{2n}(r)dr]\rightarrow 0 \ q.s., \ as \
n\rightarrow\infty,$$ we have
$E_{P_{\eta^n}}(\int_s^t(h_r-\lambda)dr|{\cal F}_s)\rightarrow0$,
$P_\eta$-a.s.,  as $n\rightarrow\infty.$

{\bf Step 1.2.3.} For any $\xi\in L^1_G(\Omega_t)$,
$E_{P_{\eta^n}}(\xi|{\cal F}_s)\rightarrow E_{P_{\eta^*}}(\xi|{\cal
F}_s)$, $P_\eta$-a.s., as $n\rightarrow\infty.$

In fact, for $\xi=\varphi(B_{s_1}-B_{s_0},\cdot\cdot\cdot,
B_{s_m}-B_{s_{m-1}})\in L_{ip}(\Omega_t)$, the conclusion is
obvious. For general $\xi\in L^1_G(\Omega_t)$, there exists a
sequence $\{\xi^m\}\subset L_{ip}(\Omega_t)$ such that
$\hat{E}[|\xi^m-\xi|]=\hat{E}[\hat{E}_s(|\xi^m-\xi|)]\rightarrow0$.
So we can assume $\hat{E}_s(|\xi^m-\xi|)\rightarrow0$ q.s..

Then, $P_\eta$-a.s., we have
\begin {eqnarray*} & &|E_{P_{\eta^n}}(\xi|{\cal F}_s)- E_{P_{\eta^*}}(\xi|{\cal
F}_s)|\\
&\leq& |E_{P_{\eta^n}}(\xi|{\cal F}_s)- E_{P_{\eta^n}}(\xi^m|{\cal
F}_s)|+|E_{P_{\eta^n}}(\xi^m|{\cal F}_s)- E_{P_{\eta^*}}(\xi^m|{\cal
F}_s)|\\
& &+|E_{P_{\eta^*}}(\xi^m|{\cal F}_s)- E_{P_{\eta^*}}(\xi|{\cal
F}_s)|\\
&\leq& 2\hat{E}_s(|\xi^m-\xi|)+|E_{P_{\eta^n}}(\xi^m|{\cal F}_s)-
E_{P_{\eta^*}}(\xi^m|{\cal F}_s)|.
\end {eqnarray*} First letting $n\rightarrow\infty$, then letting
$m\rightarrow\infty$, we have $E_{P_{\eta^n}}(\xi|{\cal
F}_s)\rightarrow E_{P_{\eta^*}}(\xi|{\cal F}_s)$, $P_\eta$-a.s.. So
combining Step 2.3 and Step 2.4, we have
\begin {eqnarray} \label {eqns8}E_{P_{\eta^*}}(\int_s^th_rdr|{\cal F}_s)=\lambda(t-s) \ P_\eta-a.s..
\end {eqnarray}

{\bf Step 1.2.4.} For $0\leq s<t\leq T$, $\eta\in L^2_{F^0}$
$\sigma\in[\underline{\sigma}, \overline{\sigma}]$, set
$\eta^\sigma=\sigma$ on $[s,t]$ and $\eta^\sigma=\eta$ on $[s,t]^c$.
We have
$$E_{P_{\eta^\sigma}}(\int_s^th_rdr|{\cal F}_{s})=\lambda(t-s) \
P_{\eta^\sigma}-a.s..$$

In fact,  Step 1.2.2-Step 1.2.3 proved the following fact: If (\ref
{eqns5}), (\ref {eqns6}) hold for some $\sigma,
\sigma'\in[\underline{\sigma}, \overline{\sigma}]$, then (\ref
{eqns8}) holds for $\sqrt{\frac{\sigma^2+\sigma'^2}{2}}$. So by
repeating the Step 1.2.2-Step 1.2.3, we get the  desired result.

{\bf Step 1.2.5.} For  any simple process $\eta\in L^2_{F^0}$,
$E_{P_\eta}(K_T)=\lambda T$.

Let  $\eta_r=\sum_{i=0}^{m-1}\eta_{t_i}1_{]t_i, t_{i+1}]}(r)\in
L^2_{F^0}$ with $\eta_{t_i}=\sum_{j=1}^{n_i}a^i_j1_{A^i_j}$ an
${\cal F}^0_{t_i}$ measurable simple function, where  $\{t_0,
\cdot\cdot\cdot, t_m\}$ is a given partition of $[0, T]$. Set
$X_t=\int_0^t\eta_rdW_r$. Let $F^X=\{{\cal F}^X_t\}$ be the
filtration generated by $X$.

Fix $0\leq i<m$. Set $\eta^{j,\varepsilon}_s=\eta_s1_{[0,
t_i+\varepsilon]}(s)+a^i_j1_{]t_i+\varepsilon, T]}(s)$ and
$X^{j,\varepsilon}_t=\int_0^t\eta^{j,\varepsilon}_s dW_s$ for
$\varepsilon>0$ small enough. Let $F^{X^{j,\varepsilon}}=\{{\cal
F}_t^{X^{j,\varepsilon}}\}$ be the  filtration generated by
$X^{j,\varepsilon}$. Then
\begin {eqnarray*}E_{P_\eta}(\int_{t_i+\varepsilon}^{t_{i+1}}h_rdr)=E_{P^0}(\int_{t_i+\varepsilon}^{t_{i+1}}h_r\circ X dr)
=E_{P^0}[E_{P^0}(\int_{t_i+\varepsilon}^{t_{i+1}}h_r\circ X dr|{\cal
F}^X_{t_i+\varepsilon})].
\end {eqnarray*}

Since $A^i_j\in{\cal F}^X_{t_i+\varepsilon}={\cal
F}^{X^{j,\varepsilon}}_{t_i+\varepsilon}$ and
$X_t=\sum_{j=0}^{n_i}X^{j,\varepsilon}_t 1_{A^i_j}$ on $[0,
t_{i+1}]$, we have
\begin {eqnarray*}
& &E_{P^0}(\int_{t_i+\varepsilon}^{t_{i+1}}h_r\circ X dr|{\cal
F}^X_{t_i+\varepsilon})\\
&=& \sum_{j=1}^{n_i}E_{P^0}(1_{A^i_j}\int_{t_i}^{t_{i+1}}h_r\circ
X^{j,\varepsilon} dr|{\cal F}^X_{t_i+\varepsilon})\\
&=& \sum_{j=1}^{n_i}1_{A^i_j}E_{P^0}(\int_{t_i}^{t_{i+1}}h_r\circ
X^{j,\varepsilon} dr|{\cal
F}^{X^{j,\varepsilon}}_{t_i+\varepsilon}).
\end {eqnarray*}
Noting  that  $$E_{P^0}(\int_{t_i+\varepsilon}^{t_{i+1}}h_r\circ
X^{j,\varepsilon}dr|{\cal
F}^{X^{j,\varepsilon}}_{t_i+\varepsilon})=E_{P_{\eta^{j,\varepsilon}}}(\int_{t_i+\varepsilon}^{t_{i+1}}h_rdr|{\cal
F}_{t_i+\varepsilon})\circ
X^{j,\varepsilon}=\lambda(t_{i+1}-t_i-\varepsilon) \ P^0-a.s.,$$ by
Step 1.2.4, we have
$E_{P_\eta}(\int_{t_i}^{t_{i+1}}h_rdr)=\lambda(t_{i+1}-t_i)$ and
$E_{P_\eta}(K_T)=\lambda T$.

{\bf Step 2.} $h\equiv\lambda$.

 Let $M_t=\int_0^th_rd\langle B\rangle_s-\int_0^t2G(h_s)ds$
and $N_t=\int_0^th_sd\langle B\rangle_s-\overline{\mu}t$. As is
mentioned in the introduction, [P07] proved that $\{M_t\}$ is a
$G$-martingale. Since $\{\int_0^th_sd\langle B\rangle_s\}$ is a
process with stationary and independent increments w.r.t. the
filtration, We know that  $\{N_t\}$ is also a $G$-martingale. Let
$L_t=\hat{E}_t(\overline{\mu}T-\overline{\sigma}^2K_T)$. Then
$\{L_t\}$ is a symmetric $G$-martingale since $K_T$ is symmetric. By
the symmetry of $\{L_t\}$ we have
$$M_t=\hat{E}_t(M_T)= \hat{E}_t(L_T+N_T)=L_t+N_t.$$
By  uniqueness of  the $G$-martingale decomposition theorem, we get
$L\equiv0$ and $h\equiv\lambda$. $\Box$

\section{Characterization of the $G$-Brownian motion }
A version of the martingale characterization for the $G$-Brownian
motion was given in [XZ09], where only symmetric $G$-martingales
with Markovian property were considered. Here we shall present a
martingale characterization in a quite different form, which is a
natural but nontrivial generalization of the classical case in a
probability space.

\noindent {\bf Theorem 4.1}(Martingale characterization of the
$G$-Brownian motion)

Let $\{M_t\}$ be a symmetric $G$-martingale with $M_T\in
L^\alpha_G(\Omega_T)$ for some $\alpha>2$ and $\{\langle
M\rangle_t\}$ a process with stationary and independent increments
(w.r.t. the filtration). Then $\{M_t\}$ is a $G$-Brownian motion;

Let $\{M_t\}$  be a $G$-Brownian motion on $(\Omega_T,
L^1_G(\Omega_T), \hat{E})$. Then there exists a positive constant
$c$ such that  $\langle M\rangle_t=c\langle B\rangle_t$.

{\bf Proof.} By  Corollary 5.2 in [Song11b], there exists $h\in
M^\alpha_G(0,T)$ such that $M_t=\int_0^th_sdB_s$. So $\langle
M\rangle_t=\int_0^t h_s^2d\langle B\rangle_s$. By Theorem 3.6, there
exists some constant $c\geq0$ such that $h^2\equiv c$. Thus by
Theorem 2.12 and Remark 2.13, $\{M_t\}$ is a $G$-Brownian motion
with $M_t$ distributed as $N(0, [c^2\underline{\sigma}^2t,
c\overline{\sigma}t])$.

On the other hand, if $\{M_t\}$  is a $G$-Brownian motion on
$(\Omega_T, L^1_G(\Omega_T)$, then $\{M_t\}$ is a symmetric
$G$-martingale. By the above arguments, we have $\langle
M\rangle_t=c\langle B\rangle_t$ for some positive constant $c$.
$\Box$

Let $${\cal H}=\{a| \ a(t)=\Sigma_{k=0}^{n-1}a_{t_k}1_{]t_k,
t_{k+1}]}(t), n\in N, 0=t_0<t_1<\cdot\cdot\cdot<t_n=T\}$$ and  $
H=\{a\in {\cal H}|\lambda[a=0]=0\}$, where $\lambda$ is the Lebesgue
measure.

\noindent {\bf Lemma 4.2} Let $\{L_t\}$ be a process with absolutely
continuous paths. Assume that there exist real numbers
$\underline{c}\leq \overline{c}$ such that  $\underline{c}(t-s)\leq
L_t-L_s\leq\overline{c}(t-s)$ for any $s<t$. Let
$C(a)=\overline{c}a^+-\underline{c}a^-$ for any $a\in R$. If
$$\hat{E}(\int_0^Ta(s)dL_s)=\int_0^TC(a(s))ds, \ \textmd{for all} \ a\in {\cal H},$$ we have that $\{L_t\}$ is a process with stationary and independent
increments such that $\underline{c}t=-\hat{E}(-L_t)\leq\hat{E}(L_t)=
\overline{c}t$, i.e., its distribution is determined by
$\underline{c}, \overline{c}$.

{\bf Proof.}  It suffices to prove the Lemma for the case
$\underline{c}< \overline{c}$. For any $a\in H$, let
$$\theta^a_s=\overline{c}1_{[a(s)\geq0]}+\underline{c}1_{[a(s)<0]}.$$
By assumption,
$$\hat{E}(\int_0^Ta(s)dL_s)=\int_0^Ta(s)\theta^a_sds.$$
On the other hand, by Theorem 2.12, there exists some weak compact
subset ${\cal P}\subset{\cal M}_1(\Omega_T)$ such that
$$\hat{E}(\xi)=\max_{P\in{\cal
P}}E_P(\xi), \ \textmd{for all} \ \xi\in L^1_G(\Omega_T),$$ which
means that there exists $P_a\in{\cal P}$ such that
$$E_{P_a}(\int_0^Ta(s)dL_s)=\int_0^Ta(s)\theta^a_sds.$$ By the
assumption for $\{L_t\}$, we have $P_a\{L_t=\int_0^t\theta^a_sds,
 \ \textmd{for all} \ t\in[0,T]\}=1$. From this we have
$$\hat{E}[\varphi(L_{t_1}-L_{t_0},\cdot\cdot\cdot,
L_{t_n}-L_{t_{n-1}})]\geq\varphi(\int_{t_0}^{t_1}\theta^a_sds,\cdot\cdot\cdot,
\int_{t_{n-1}}^{t_n}\theta^a_sds)$$ for any $\varphi\in C_b(R^n)$
and $n\in N$. Consequently,
\begin {eqnarray*}& &\hat{E}[\varphi(L_{t_1}-L_{t_0},\cdot\cdot\cdot, L_{t_n}-L_{t_{n-1}})]\\
                  &\geq&\sup_{a\in H}\varphi(\int_{t_0}^{t_1}\theta^a_sds,\cdot\cdot\cdot,\int_{t_{n-1}}^{t_n}\theta^a_sds)\\
                  &=&\sup_{c_1,\cdot\cdot\cdot,c_n\in[\underline{c},\overline{c}]}\varphi(c_1(t_1-t_0),\cdot\cdot\cdot,c_n(t_n-t_{n-1})).
\end {eqnarray*} The converse inequality is obvious. Thus $\{L_t\}$ is a process with stationary and independent
increments such that $\underline{c}t=-\hat{E}(-L_t)\leq\hat{E}(L_t)=
\overline{c}t$. $\Box$

\noindent {\bf Lemma 4.3} Let $\{L_t\}$ be a $G$-martingale with
finite variation and $L_T\in L^\beta_G(\Omega_T)$ for some
$\beta>1$. Then $\{L_t\}$ is non-increasing.  Particularly, $L_t\leq
L_0=\hat{E}(L_T).$

{\bf Proof.} By Theorem 4.5 in [Song11a], we know $\{L_t\}$ has the
following decomposition
$$L_t=\hat{E}(L_T)+M_t+K_t,$$ where $\{M_t\}$ is a symmetric
$G$-martingale and $\{K_t\}$ is a non-positive, non-increasing
$G$-martingale. Since both $\{L_t\}$ and $\{K_t\}$ are processes
with finite variation, we get $M_t\equiv0$. Therefore, we have
$L_t=\hat{E}(L_T)+K_t\leq\hat{E}(L_T)=L_0.$ $\Box$

\noindent {\bf Theorem 4.4} Let $\{X_t\}$ be a generalized
$G$-Brownian motion with zero mean. Then we have the following
decomposition:
$$X_t=M_t+L_t,$$ where $\{M_t\}$ is a symmetric $G$-Brownian motion,
and $\{L_t\}$ is a non-positive, non-increasing $G$-martingale with
stationary and independent increments.

{\bf Proof.} Clearly $\{X_t\}$ is a  $G$-martingale. By Theorem 4.5
in [Song11a], we have the following decomposition $$X_t=M_t+L_t,$$
where $\{M_t\}$ is a symmetric $G$-martingale, and $\{L_t\}$ is a
non-positive, non-increasing $G$-martingale. Noting that $X_t\in
L^3_G(\Omega_T)$ from the definition of generalized $G$-Brownian
motion, we know that $M_t, L_t\in L^\beta_G(\Omega_T)$ for any
$1\leq\beta<3$ by Theorem 4.5 in [Song11a].

In the sequel, we  first prove that $\{L_t\}$ is a process with
stationary  and independent increments. Noting that
$\hat{E}(-L_t)=\hat{E}(-X_t)=ct$ for some positive constant $c$
since $\{X_t\}$ is a process with stationary and independent
increments, we claim that $-L_t-ct$ is a $G$-martingale. To prove
this, it suffices to show that for any $t>s$,
$\hat{E}_s[-(L_t-L_s)]=c(t-s).$ In fact, since $\{M_t\}$ is a
symmetric $G$-martingale, we have
 $$\hat{E}_s[-(L_t-L_s)]=\hat{E}_s[-(X_t-M_t-X_s+M_s)]=\hat{E}_s[-(X_t-X_s)].$$

Noting that $\{X_t\}$ is a process with independent
increments(w.r.t. the filtration),
$$\hat{E}_s[-(X_t-X_s)]=\hat{E}[-(X_t-X_s)]=c(t-s).$$

Combining this with Lemma 4.3, we have $-(L_t-L_s)-c(t-s)\leq0$ for
any $s<t$. On the other hand, for any $a\in {\cal H}$, noting that
$\{M_t\}$ is a symmetric $G$-martingale, we have
$$\hat{E}[\int_0^Ta(s)dL_s]=\hat{E}[\int_0^Ta(s)dX_s]=\hat{E}[\Sigma_{k=0}^{n-1}a_{t_k}(X_{t_{k+1}}-X_{t_k})].$$
Since $\{X_t\}$ is a process with stationary, independent
increments, we have

\begin {eqnarray*}& &\hat{E}[\int_0^Ta(s)dL_s]\\
                  &=&\Sigma_{k=0}^{n-1}\hat{E}[a_{t_k}(X_{t_{k+1}}-X_{t_k})]\\
                  &=&\Sigma_{k=0}^{n-1}ca^-_{t_k}(t_{k+1}-t_k)\\
                  &=&\int_0^Tca^-(s)ds=\int_0^TC(a(s))ds,
\end {eqnarray*}
where $C(a(s))$ is defined as in Lemma 4.2 with $\overline{c}=0,
\underline{c}=-c$. By Lemma 4.2, $\{L_t\}$ is a process with
stationary  and independent increments.

Now we are in a position to show that $\{M_t\}$ is a (symmetric)
$G$-Brownian motion. To this end, by Theorem 4.1, it suffices to
prove that $\{\langle M\rangle_t\}$ is a process with stationary and
independent increments (w.r.t. the filtration). For $n\in N$, let
$$X^n_t=\sum_{k=0}^{2^n-1}X_{\frac{kT}{2^n}}1_{]\frac{kT}{2^n},
\frac{(k+1)T}{2^n}]}(t)$$ and
$$\Omega^n_t(X)=\sum_{k=0}^{2^n-1}(X_{\frac{(k+1)t}{2^n}}-X_{\frac{kt}{2^n}})^2.$$

Observing that $\Omega^n_t(X)=X_t^2-2\int_0^tX^n_sdX_s$, we have
\begin {eqnarray*}
& &|\Omega^n_t(X)-\Omega^{m+n}_t(X)|\\
&\leq& 2(|\int_0^t(X^n_s-X^{m+n}_s)dM_s|+|\int_0^t(X^n_s-X^{m+n}_s)dL_s|)\\
&=&2(|I|+|II|).
\end {eqnarray*} for any $n,m\in N.$  It's easy to check that
$$\hat{E}(|II|)\leq c\int_0^t\hat{E}(|X^n_s-X^{m+n}_s|)ds\rightarrow0 \ \textmd{as} \ m, n\rightarrow\infty.$$
Noting that
\begin {eqnarray*}I&=&\sum_{i=0}^{2^n-1}\sum_{j=0}^{2^m-1}(X_{\frac{it}{2^n}+\frac{jt}{2^{n+m}}}-X_{\frac{it}{2^n}})
(M_{\frac{it}{2^n}+\frac{(j+1)t}{2^{n+m}}}-M_{\frac{it}{2^n}+\frac{jt}{2^{n+m}}})\\
&=& \sum_{i=0}^{2^n-1}\sum_{j=0}^{2^m-1}I_i^j,
\end {eqnarray*} we get
$$\hat{E}(I^2)\leq\sum_{i=0}^{2^n-1}\sum_{j=0}^{2^m-1}\hat{E}[(I_i^j)^2].$$

Let's estimate the expectation $\hat{E}[(I_i^j)^2]$:
\begin {eqnarray*}& &\hat{E}[(I_i^j)^2]\\
&=&\hat{E}[(X_{\frac{it}{2^n}+\frac{jt}{2^{n+m}}}-X_{\frac{it}{2^n}})^2
(M_{\frac{it}{2^n}+\frac{(j+1)t}{2^{n+m}}}-M_{\frac{it}{2^n}+\frac{jt}{2^{n+m}}})^2]\\
&\leq&2
\hat{E}[(X_{\frac{it}{2^n}+\frac{jt}{2^{n+m}}}-X_{\frac{it}{2^n}})^2
\{(X_{\frac{it}{2^n}+\frac{(j+1)t}{2^{n+m}}}-X_{\frac{it}{2^n}+\frac{jt}{2^{n+m}}})^2+
(L_{\frac{it}{2^n}+\frac{(j+1)t}{2^{n+m}}}-L_{\frac{it}{2^n}+\frac{jt}{2^{n+m}}})^2\}]
\end {eqnarray*}
Noting that $-c(t-s)\leq L_t-L_s\leq0$, we have
$$\hat{E}[(I_i^j)^2]\leq
\hat{E}[(X_{\frac{it}{2^n}+\frac{jt}{2^{n+m}}}-X_{\frac{it}{2^n}})^2
\{(X_{\frac{it}{2^n}+\frac{(j+1)t}{2^{n+m}}}-X_{\frac{it}{2^n}+\frac{jt}{2^{n+m}}})^2+
c^2\frac{t^2}{2^{2(n+m)}}\}].$$

By (\ref {eqns7}), $\hat{E}[(X_t-X_s)^2]\leq C_1|t-s|$ for some
constant $C_1$. From the condition of independent increments of $X$,
we have  $\hat{E}[(I_i^j)^2]\leq C\frac{j}{2^{2(n+m)}}$ for some
constant $C$, hence that $\hat{E}(I^2)\rightarrow0$, and finally
that $\hat{E}(|\Omega^n_t(X)-\Omega^{m+n}_t(X)|)\rightarrow0$ as
$m,n\rightarrow\infty$. Then $$\langle
X\rangle_t:=\lim_{L^1_G(\Omega_T), n\rightarrow\infty}\Omega^n_t$$
is a process with stationary and independent increments (w.r.t. the
 filtration). Noting that $\langle M\rangle_t=\langle X\rangle_t$, $\langle
M\rangle_t$ is also a process with stationary and independent
increments (w.r.t. the
 filtration). $\Box$
\section{$G$-martingales with finite variation}

\noindent {\bf Proposition 5.1} Let $\eta\in M^1_G(0,T)$ with
$|\eta|\equiv c$ for some constant $c$. Then
\begin {eqnarray}\label{eqny1} K_t:=\int_0^t\eta_sd\langle
B\rangle_s-\int_0^t2G(\eta_s)ds
\end {eqnarray}
is a process with stationary and independent increments. Moreover,
for fixed $c$, all processes in the above form have the same
distribution.

{\bf Proof.} Since
$-c(\overline{\sigma}^2-\underline{\sigma}^2)(t-s)\leq K_t-K_s\leq0$
for any $s<t$, by Lemma 4.2, it suffices to prove that for any
$a\in{\cal H}$
$$\hat{E}(\int_0^Ta_sdK_s)=\int_0^TC(a_s)ds,$$ where $C(a_s)$ is defined as in Lemma 4.2 with
$\overline{c}=0,
\underline{c}=-c(\overline{\sigma}^2-\underline{\sigma}^2)$. In
fact, noting that $$\int_0^Ta_sdK_s\leq
\int_0^T2G(a_s\eta_s)ds-\int_0^T2a_sG(\eta_s)ds=\int_0^TC(a_s)ds,$$
we have
$$\hat{E}(\int_0^Ta_sdK_s)\leq\int_0^TC(a_s)ds.$$ On the other
hand, we have
$$\hat{E}(\int_0^Ta_sdK_s)\geq-\hat{E}\{-[\int_0^T2G(a_s\eta_s)ds-\int_0^T2a_sG(\eta_s)ds]\}=\int_0^TC(a_s)ds.$$ So $\{K_t\}$
is a process with stationary and independent increments and its
distribution is determined by $c$. $\Box$

Just like the conjecture by Shige Peng for the representation of
 $G$-martingales with finite variation, we guess that any
$G$-martingale with  stationary, independent increments and finite
variation should have the form of (\ref {eqny1}). At the end we
present a characterization for $G$-martingales with finite
variation.

\noindent {\bf Proposition 5.2} Let $\{M_t\}$ be a $G$-martingale
with $M_T\in L^\beta_G(\Omega_T)$ for some $\beta>1$. Then $\{M_t\}$
is a $G$-martingale with finite variation if and only if
$\{f(M_t)\}$ is  a $G$-martingale for any non-decreasing $f\in
C_{b,lip}(R)$.

{\bf Proof.} Necessity. Assume $\{M_t\}$ is a $G$-martingale with
finite variation. By Lemma 4.3, we know that $\{M_t\}$ is
non-increasing. By Theorem 5.4 in [Song11b], there exists a sequence
$\{\eta^n_t\}\subset H^0_G(0,T)$ such that
$$\hat{E}[\sup_{t\in[0,T]}|M_t-L_t(\eta^n)|^\beta]\rightarrow0$$ as $n$ goes
to infinity, where $L_t(\eta^n)=\int_0^t\eta^n_sd\langle
B\rangle_s-\int_0^t2G(\eta^n_s)ds$. It suffices to prove that for
any $\eta\in H^0_G(0,T)$ and non-decreasing $f\in C_b^2(R)$,
$f(L_t(\eta))$ is a $G$-martingale. In fact,
\begin {eqnarray*}f(L_t(\eta))&=&f(L_0)+\int_0^tf'(L_s(\eta))dL_s(\eta)\\
                              &=&f(L_0)+\int_0^t f'(L_s(\eta))\eta_sd\langle B\rangle_s-\int_0^t2
f'(L_s(\eta))G(\eta_s)ds.
\end {eqnarray*} Since $f'(L_s(\eta))\geq0$ and $f'(L_s(\eta))\eta_s\in
M^1_G(0,T)$, we conclude that
$$f(L_t(\eta))=f(L_0)+L_t(f'(L(\eta))\eta)$$ is a $G$-martingale.

Sufficiency. Assume $\{f(M_t)\}$ is  a $G$-martingale for any
non-decreasing $f\in C_{b,lip}(R)$. Let $X_t:=\arctan M_t$. Then
$\{X_t\}$ is a bounded $G$-martingale and $\{f(X_t)\}$ is  a
$G$-martingale for any non-decreasing $f\in C_{b,lip}(R)$. By
Theorem 4.5 in [Song11a], we know $\{X_t\}$ has the following
decomposition
$$X_t=\hat{E}(X_T)+N_t+K_t,$$ where $\{N_t\}$ is a symmetric
$G$-martingale and $\{K_t\}$ is a non-positive, non-increasing
$G$-martingale. Then by It\^{o}'s formula $$e^{\alpha X_t}=e^{\alpha
X_0}+\alpha\int_0^t e^{\alpha X_s}dX_s+\frac{\alpha^2}{2}\int_0^t
e^{\alpha X_s}d\langle N\rangle_s.$$ For any $\alpha>0$, by
assumption, $e^{\alpha X_t}$ is a $G$-martingale. So $L_t:=\int_0^t
e^{\alpha X_s}dK_s+\frac{\alpha}{2}\int_0^t e^{\alpha X_s}d\langle
N\rangle_s$ is a $G$-martingale  with finite variation. By Lemma
4.3, $L_t$ is non-increasing, by which we conclude that
$K_t+\frac{\alpha}{2}\langle N\rangle_t$ is non-increasing. So
$$\frac{\alpha}{2}\hat{E}(\langle N\rangle_T)\leq\hat{E}(-K_T) \ \textmd{for all} \ \alpha>0.$$
By this, we conclude that $\hat{E}(\langle N\rangle_T)=0$ and
$N_t\equiv0$. Then $X_t=\hat{E}(X_T)+K_t$ is non-increasing, and
consequently, $M_t$ is non-increasing. $\Box$

Particularly, Proposition 5.2 provides a method to convert
$G$-martingales with finite variation into bounded $G$-martingales
with finite variation.




\providecommand{\bysame}{\leavevmode\hbox
to3em{\hrulefill}\thinspace}
\providecommand{\MR}{\relax\ifhmode\unskip\space\fi MR }
\providecommand{\MRhref}[2]{%
  \href{http://www.ams.org/mathscinet-getitem?mr=#1}{#2}
} \providecommand{\href}[2]{#2}

\end{document}